\documentclass[12pt]{article}
\usepackage{amsmath}
\usepackage{amssymb}
\usepackage{url}

\usepackage[T1]{fontenc}
\usepackage[light,math]{anttor}
\usepackage{diagrams}
\diagramstyle{midshaft,PostScript=dvips,nohug}

\usepackage{theorem}
\theoremstyle{change}
\theorembodyfont{\normalfont\itshape}
\newtheorem{lemma}{Lemma.}[section]
\newtheorem{prop}[lemma]{Proposition.}
\newtheorem{thm}[lemma]{Theorem.}
\newtheorem{cor}[lemma]{Corollary.}

\newtheorem{taller}[lemma]{$\!\!$}
\newenvironment{blanko}[1]%
{\begin{taller}{\bf{#1}}\normalfont}%
{\end{taller}}

\providecommand{\qed}{\hspace*{\fill}\nolinebreak[1]\hspace*{\fill}$\Box$}

\newenvironment{proof}{\begin{list}{\em Proof. }%
{\setlength{\labelsep}{0mm}\setlength{\leftmargin}{0mm}%
\setlength{\labelwidth}{0mm}\setlength{\listparindent}{\parindent}%
\setlength{\parsep}{\parskip}\setlength{\partopsep}{0mm}}%
\item}{\hfill \qed\end{list}}

\newcommand{\Fib}{\mathbf{Fib}}
\newcommand{\BiFib}{\mathbf{BiFib}}
\newcommand{\ev}{\operatorname{ev}}
\newcommand{\id}{\operatorname{id}}
\newcommand{\dom}{\operatorname{dom}}

\newcommand{\upperstar}{\raisebox{0.7ex}[0ex][0ex]{\ensuremath{\ast}}}
\newcommand{\lowerstar}{{\raisebox{-0.6ex}[-0.5ex][0ex]{\ensuremath{\ast}}}}
\newcommand{\lowershriek}{_{\,!\,}}

\newcommand{\isopil}{\stackrel{\raisebox{0.1ex}[0ex][0ex]{\(\sim\)}}%
           {\raisebox{-0.15ex}[0.28ex]{\(\rightarrow\)}}}

\newcommand{\E}{\mathcal{E}}
\newcommand{\F}{\mathcal{F}}
\newcommand{\G}{\mathcal{G}}
\newcommand{\Set}{\mathbf{Set}}

\begin{document}
%%%%%%%%%%%%%%%%%%%%%%%%%%%%%%%%%%%%%%%%%%%%%%%%%%
\title{Local fibered right adjoints are polynomial}
\author{Anders Kock 
 \and 
 Joachim Kock\thanks{Partially supported by research grants 
  MTM2009-10359 % Nart
       and
MTM2010-20692 % Castellana
of the Spanish Ministry of Science and Innovation.}
}
% \author{Anders Kock \\ \footnotesize{ Matematisk Institut} \\[-4pt] \footnotesize{Aarhus Universitet, Denmark}
% \\[-4pt]
% \footnotesize{
%  \url{kock@imf.au.dk}} 
%  \and 
%  Joachim Kock\thanks{Partially supported by research grants 
%   MTM2009-10359 % Nart
%        and
% MTM2010-20692 % Castellana
% of the Spanish Ministry of Science and Innovation.}
% \\
% \footnotesize{Departament de Matem\`atiques}
% \\[-4pt]
% \footnotesize{Universitat Aut\`onoma de Barcelona, 
% Spain}
% \\[-4pt]
% \footnotesize{\url{kock@mat.uab.cat}}}

\date{}

\maketitle

\small

\noindent {\bf Abstract.} For any locally cartesian closed category $\E$, we prove
  that a local fibered right adjoint between slices of $\E$ is given by a
  polynomial.  The slices in question are taken in a well known fibered sense.

\normalsize

% MSC2000:
%   18D30 Fibered categories
%   18D15 Closed categories (closed monoidal and Cartesian closed categories, etc.)

% Keywords:
%   fibered categories,
%   locally cartesian closed categories,
%   polynomial functors

%%%%%%%%%%%%%%%%%%%%%%%%%%%%%%%%%%%%%%%%%%%%%%%%%%
\section{Introduction}
%%%%%%%%%%%%%%%%%%%%%%%%%%%%%%%%%%%%%%%%%%%%%%%%%%

In a locally cartesian closed category $\E$, a diagram of the form
\begin{equation}\begin{diagram}[w=4ex,tight]
I& \lTo^{s}&E& \rTo^{p}&B&\rTo^{t} &J
\end{diagram}\label{bridge}\end{equation}
gives rise to a so-called {\em polynomial functor} $\E/I \to \E/J$, namely the
composite functor
\begin{equation}\label{plainpolyfun}\begin{diagram}
\E/I & \rTo^{s\upperstar}&\E/E&\rTo^{p\lowerstar}&\E/B&\rTo^{t\lowershriek}&\E/J.
\end{diagram}\end{equation}

The use of polynomial functors as data type constructors goes back at least to
the 1980s, cf.~Manes and Arbib (1986).  Moerdijk and
Palmgren (2000) (and in a more general setting,
Gambino and Hyland (2004)) showed that initial algebras for
polynomial endofunctors are precisely the W-types of Martin-L\"of type theory. 
In a series of papers Abbott, Altenkirch, Ghani and collaborators have
  developed further the theory of polynomial functors (called container
  functors in (Abbott {\em et al.} 2003)) and their natural transformations
  as data type constructors and polymorphic functions, subsuming notions as
  shapely types, strictly positive types, and general tree types.
%   %
% In
% a series of papers, Abbott, Altenkirch, Ghani (2003--) and
% collaborators have developed further the theory of polynomial functors --- they
% call them container functors --- and their natural transformations as data type
% constructors and polymorphic functions, subsuming notions as shapely types,
% strictly positive types, and general tree types.
% %
We refer to
(Gambino-Kock 2009) for background on polynomial functors
and an extensive bibliography, including also pointers to the use of polynomial
functors in logic, combinatorics, representation theory, topology, and higher
category theory.

In (Gambino-Kock 2009), six different intrinsic characterizations of
polynomial functors are listed for the case where $\E = \Set$,
one of them being that a functor $P:\Set/I\to\Set/J$ is polynomial if and only
if it is a local right adjoint, i.e.~the slice of $P$ at the terminal object of
$\Set/I$ is a right adjoint.  Already for $\E$ a presheaf topos, this 
characterization
fails, as pointed out by Weber (2007).  His counter-example is a
significant one: the free-category monad on the category of graphs is a local
right adjoint, but it is not polynomial.

In this note we adjust the local-right-adjoint characterization so as to be
valid in every locally cartesian closed category; for this we pass to the
setting of fibered slice categories and fibered functors.  The {\em fibered
slice} $\E|I$ is the fibered category over $\E$ whose fiber over an object
$K$ is the plain slice $\E/(I\!\times\! K)$.  A diagram as \eqref{bridge} defines
also a fibered polynomial functor
\begin{diagram}\E|I & \rTo ^{s\upperstar}&\E|E&\rTo ^{
p\lowerstar}&\E|B&\rTo^{t\lowershriek}&\E|J,
\end{diagram}
whose fiber over the terminal object $1\in\E$ is the plain polynomial functor \eqref{plainpolyfun}.  Hence any
polynomial functor has a canonical extension to a fibered functor.

Our main theorem is this: {\em a fibered functor $P:\E|I \to \E|J$ is
polynomial if and only if it is a local fibered right adjoint}.

%%%%%%%%%%%%%%%%%%%%%%%%%%%%%%%%%%%%%%%%%%%%%%%%%%
\section{Fibered categories and fibered slices}
%%%%%%%%%%%%%%%%%%%%%%%%%%%%%%%%%%%%%%%%%%%%%%%%%%

This section is mainly to fix terminology and notation.
We refer to (Borceux 1994), (Johnstone 2002), and (Streicher 1999) for
further background on fibered categories. 
   
A category $\E$ is fixed throughout; for the main result, it should
be a locally cartesian closed category (lccc) with a terminal object.
For some of the considerations,  it suffices that it is a category with finite
limits, as in (Streicher 1999).
We also assume we have chosen
pullbacks once and for all, so that for each arrow $a: J\to I$
we have at our disposal two functors
$$
a\lowershriek : \E/J \rTo \E/I
\qquad
a\upperstar : \E/I \rTo \E/J,
$$
and if $\E$ is furthermore a lccc, the $a\upperstar$s have right adjoints
$$
a\lowerstar : \E/J \rTo \E/I
$$
which we also assume are chosen once and for all; so
$$
a\lowershriek \dashv a\upperstar  \dashv a\lowerstar .
$$

\begin{blanko}{Fibered categories.}
  We shall work
with categories fibered over $\E$, henceforth just called
 fibered categories.  These form a $2$-category $\Fib_\E$ whose objects are
%%??
 fibered cate\-gories, whose $1$-cells are fibered functors (i.e.~functors
 commuting with the structure functor to $\E$ and sending cartesian arrows to
 cartesian arrows), and whose $2$-cells are fibered natural transformations
 (i.e.~natural transformations whose components are vertical arrows).
  If $\F$ is a fibered category, we denote by $\F^I$ the (strict)
  fiber over an object $I\in \E$, and if $L:\G\to\F$ is a fibered
  functor, we write $L^I:\G^I \to \F^I$ for the induced functor
  between the $I$-fibers.
\end{blanko}

\begin{blanko}{Cleavage and base change.}
  A {\em cleavage} of a fibered category $\F$ is a choice of cartesian lifts:
  for each arrow $a:J \to I$ in $\E$ and object $X\in \F^I$ we choose a
  cartesian arrow over $a$ with codomain $X$.  For each $a$ this assignment
  defines a functor $a\upperstar:\F^{I}\to \F^{J}$ called {\em base change}
  along $a$: the value of $a\upperstar$ on $X\in \F^{I}$ is taken to
  be the domain of the chosen cartesian lift of $a$ with codomain $X$, and the
  value on morphisms comes about by exploiting the universal property of
  cartesian arrows.

  For a composable pair of arrows $a_{1}, a_{2}$ in $\E$, there is a canonical
  isomorphism $a_{1}\upperstar\circ a_{2}\upperstar \isopil (a_{2}\circ
  a_{1})\upperstar$, likewise derived from the universal property of cartesian
  arrows.

  For a fibered functor $L: \G \to \F$, the functors induced on the fibers
  commute with base change functors up to canonical isomorphisms: for $a:J\to I$
  in $\E$, the universal property of the chosen cartesian lifts of $a$ in $\F$
  gives rise to 2-cells
  \begin{equation}\label{upperstar2cell}
  \begin{diagram}[w=5ex,h=5ex,tight]
  \G^I & \rTo^{L^I} & \F^I  \\
  \dTo<{a\upperstar} & L^a\!\raisebox{-1ex}{\ensuremath{\nearrow}} & \dTo>{a\upperstar} \\
  \G^{J} & \rTo_{L^{J}} & \F^{J}
  \end{diagram}
  \end{equation}
   The fact that the $L^{a}$s are invertible comes from the fact that $L$
   preserves cartesian arrows.  Again, because the isomorphisms come about from
   universal properties, coherence conditions can be deduced.
\end{blanko}

\begin{blanko}{Fibered adjunctions.}
 An adjunction between two fibered functors in opposite directions is given by
 two fibered natural transformations $\eta: \id \Rightarrow R\circ L$ and
 $\varepsilon: L \circ R \Rightarrow \id$ satisfying the usual triangle
 identities.  In other words, it is an adjunction in the $2$-category
 $\Fib_\E$.  It is clear that a fibered adjunction induces a fiber-wise adjunction
 each fiber.  Conversely (cf.~e.g.~(Borceux 1994) 8.4.2), if $R : \F\to \G$
 is a fibered functor and for each $I\in \E$ there is a left adjoint $L^I
 \dashv R^I$, then these $L^I$ assemble into a fibered left adjoint if for
 every arrow $a: J \to I$ in $\E$, the mate of the canonical invertible
 $2$-cell
 \begin{diagram}[w=5.5ex,h=5.5ex,tight]
 \F^I & \rTo^{R^I} & \G^I  \\
 \dTo<{a\upperstar} & \raisebox{0.5ex}{\ensuremath{(R^a)^{-1}}}
 \!\!\!\!\!\!\raisebox{-1.7ex}{\ensuremath{\swarrow}} & \dTo>{a\upperstar} \\
 \F^{J} & \rTo_{R^{J}} & \G^{J} ,
 \end{diagram}
 namely
 \begin{diagram}[w=5ex,h=5ex,tight]
 \F^I & \lTo^{L^I} & \G^I  \\
 \dTo<{a\upperstar} & \nwarrow & \dTo>{a\upperstar}  \\
 \F^{J} & \lTo_{L^{J}} & \G^{J}
 \end{diagram}
 is again invertible, and hence turns the family $L^I$ into a fibered functor.
\end{blanko}

\begin{blanko}{Bifibrations and $\E$-indexed sums.}\label{bifibrations}
    A fibered category $\F$ is called {\em bifibered} if the structure functor
    $\F\to\E$ is also an opfibration, i.e.~has all opcartesian lifts.  We then
    assume the choice of an opcartesian lift for each arrow $a:J\to I$ and each
    object $T\in \F^J$.  This defines cobase-change functors, which we denote by
    lowershriek: $a\lowershriek : \F^J \to \F^I$.  Cobase change is left adjoint
    to base change: $a\lowershriek \dashv a\upperstar $.  Indeed, for $a: J \to
    I$ in $\E$ we have natural bijections $$ \F^J(T, a\upperstar X) \cong
    \F^a(T,X) \cong \F^I(a\lowershriek T, X) $$ according to the universal
    properties of cartesian and opcartesian arrows.  (Here $\F^a(T,X)$ denotes
    the set of arrows $T\to X$ lying over $a$.)

 A fibered functor is called {\em bifibered} if it preserves also opcartesian arrows.
 In terms of cobase-change functors, we can say that a fibered functor $L:
 \G\to\F$ is bifibered if for every arrow $a:J\to I$ in $\E$, the mate of the
 compatibility-with-base-change square \eqref{upperstar2cell}:
   \begin{equation}\label{lowershriek2cell}
 \begin{diagram}[w=5ex,h=5ex,tight]
 \G^I & \rTo^{L^I} & \F^I  \\
 \uTo<{a\lowershriek} & \raisebox{-1ex}{\ensuremath{\nwarrow}}\overline{L}{}^a
 & \uTo>{a\lowershriek} \\
 \G^{J} & \rTo_{L^{J}} & \F^{J}
 \end{diagram}
 \end{equation}
 is invertible.
 
  Bifibered categories are mostly interesting when they further satisfy
  the {\em Beck-Chevalley condition}.
% \end{blanko}
% 
% \begin{blanko}{$\E$-indexed sums.}
\label{E-indexedsums}
  In terms of chosen cartesian and
  opcartesian lifts this condition says that for every pullback square in $\E$
  \begin{diagram}[w=4ex,h=4ex,tight]
  \cdot\SEpbk & \rTo^b & \cdot \\
  \dTo<u & & \dTo>v \\
  \cdot & \rTo_a & \cdot
  \end{diagram}
  the fibered natural transformation
  $$
  u\lowershriek b\upperstar \ \Rightarrow \ a\upperstar
  v\lowershriek
  $$
  is invertible.
 
  A fibered category $\F$ is said to have {\em $\E$-indexed sums} when $\F\to\E$
  is bifibered, and the Beck-Chevalley condition holds. 
 %%?
 We denote by $\BiFib_\E$ the full subcategory of bifibered 
categories/bifibered functors 
thus determined.
%   between categories with $\E$-indexed sums is called {\em $\E$-linear}.
\end{blanko}

\begin{prop}\label{fibleftpreslowershriek}
  A fibered left adjoint between bifibered categories 
preserves cobase change (in particular, preserves $\E$-indexed sums if 
the categories have such).
\end{prop}
% \begin{prop}
%   A fibered left adjoint between fibered categories 
% preserves cobase change (in particular, preserves $\E$-indexed sums if 
% the categories have such).
% \end{prop}
\begin{proof}
Take left adjoints of all the arrows in the
base-change compatibility square \eqref{upperstar2cell} for the right adjoint.
\end{proof}
% (In fact, $\E$-indexed sums are examples of fibered colimits, and fibered left
% adjoints preserve all fibered colimits.)

\begin{blanko}{Fibered slices.}
The {\em fibered slice} $\E|I$ is the category whose objects are spans
\begin{diagram}[w=4ex,tight]
I & \lTo^p & M & \rTo^q & K ,
\end{diagram}
and whose morphisms are diagrams
\begin{diagram}[w=5ex,h=4ex,tight,nohug]
& & M' & \rTo^{q'} & K' \\
I &\ldTo(2,1)^{p'} & \dTo>v && \dTo>w \\
& \luTo(2,1)_p & M & \rTo_{q} & K .
\end{diagram}
 The structural functor $\E|I \to \E$ that returns the right-most 
object (resp.\ arrow) is
a bifibration.  The cartesian arrows are the diagrams for which the square is
a pullback, while the opcartesian arrows are those for which $v$ is 
invertible, as is easy to verify.
The vertical arrows are those for which $w$ is an identity arrow.

More conceptually, the fibered slice is obtained from the plain slice as
the following pullback, and the structural functor as the left-hand vertical
composite:
\begin{diagram}[w=5ex,h=5ex,tight]
\E|I \SEpbk & \rTo  & \E/I  \\
\dTo  &    & \dTo>{\dom}  \\
\operatorname{Ar}(\E)  & \rTo^d  & \E \\
\dTo<c \\
\E .
\end{diagram}
(Here $\operatorname{Ar}(\E)$ is the category of arrows in $\E$, and $d$ and $c$
are the domain and codomain fibrations, respectively.)  This is to say that the
fibered slice is the so-called family fibration of the fibration
$\dom:\E/I\to\E$, cf.~(Streicher 1999) 6.2 
%%? 
or (Johnstone 2002) Prop.\ 1.4.16.

For the $K$-fiber of $\E|I$  we have the canonical identification
$$
(\E|I)^K \cong \E/(I\!\times\! K) .
$$
In particular, the plain slice $\E /I$ sits
inside the fibered slice $\E|I$ as the fiber over
the terminal object $1\in \E$:
$$
 (\E|I)^{1} \cong \E/I .
$$
Note also that we have $\E|1 \cong \operatorname{Ar}(\E)$.

In the $I$-fiber of $\E|I$ we have the canonical object given by the
identity span
\begin{diagram}[w=4ex,tight]
I & \lTo^\id & I & \rTo^\id & I.
\end{diagram}
which, in view of the interpretation in terms of plain slices, we denote by
$\delta$ (for ``diagonal'').

We note that $\E|I$ has $\E$-indexed sums: for $a:K'\to K$ in $\E$, the
base-change functor $a\upperstar : (\E|I)^K \to (\E|I)^{K'}$ is identified
with $(\id_I\times a)\upperstar : \E/(I\!\times\!  K) \to \E/(I\!\times\!
K')$ which has left adjoint $(\id_I\times a)\lowershriek$, and the Beck-Chevalley
condition follows from the case of plain slices.
\end{blanko}

\begin{blanko}{Polynomial functors --- fibered version.}\label{polyparagraph}
We now assume that $\E$ is a lccc. Each arrow $f:J \to I$ in $\E$ induces
fibered functors 
% (the third provided $\E$ is a lccc)
$$
f\lowershriek : \E|J \rTo \E|I
\qquad
f\upperstar : \E|I \rTo \E|J
\qquad
f\lowerstar : \E|J \rTo \E|I
$$
and fibered adjunctions
$$
f\lowershriek \dashv f\upperstar \dashv f\lowerstar .
$$
These extend the basic functors on plain slices: for example if $f\upperstar :
\E/I \to \E/J$ is the plain pullback, then the $K$-fiber of the
fibered pullback functor $f\upperstar$ is
$$
(f \!\times\! \id_K) \upperstar: \E/(I\!\times\! K) \to \E/(J\!\times\! K) .
$$

A fibered functor of the form
$$\begin{diagram}\E|I & \rTo ^{s\upperstar}&\E|E&\rTo ^{
p\lowerstar}&\E|B&\rTo^{t\lowershriek}&\E|J
\end{diagram}$$
for a diagram in $\E$ 
\begin{diagram}[w=4ex,tight]
I& \lTo ^{s}&E& \rTo ^{p}&B&\rTo ^{t} &J
\end{diagram}
is called a (fibered) {\em polynomial functor}.
\end{blanko}

%%%%%%%%%%%%%%%%%%%%%%%%%%%%%%%%%%%%%%%%%%%%%%%%%%
\section{Fibered left adjoints}
%%%%%%%%%%%%%%%%%%%%%%%%%%%%%%%%%%%%%%%%%%%%%%%%%%
\label{Sec:leftadj}

In this section, the base category $\E$ is just assumed to have finite limits.

\bigskip

Recall that $\delta\in\E|I$ denotes the identity span 
\begin{diagram}[inline,w=4ex,tight]
I & \lTo^\id & I & \rTo^\id & I.
\end{diagram}

\begin{blanko}{Main lemma.}
Let $\F$ be a fibered category with $\E$-indexed sums.
 For each $I \in \E$, the functor
 \begin{eqnarray*}
  \ev_\delta: \BiFib_\E(\E|I, \F) & \longrightarrow & \F^I \\
   L & \longmapsto & L(\delta)
 \end{eqnarray*}
 is an equivalence of categories.  A pseudo-inverse is given on objects by
 \begin{eqnarray*}
    h: \F^I   & \longrightarrow & \BiFib_\E(\E|I,\F)  \\
   X & \longmapsto & [\langle p,q\rangle \mapsto q\lowershriek p\upperstar X]
 \end{eqnarray*}
 where $\langle p,q\rangle$ denotes a span $
%   \begin{equation}\label{spann}
    \begin{diagram}[w=4ex,tight]
    I & \lTo^p  & M & \rTo^q & K .
  \end{diagram}
% \end{equation}
$
\end{blanko}
% \begin{blanko}{Main lemma.}
% Let $\F$ be a fibered category with $\E$-indexed sums.
%  For each $I \in \E$, the functor
%  \begin{eqnarray*}
%   \ev_\delta: \BiFib_\E(\E|I, \F) & \longrightarrow & \F^I \\
%    L & \longmapsto & L(\delta)
%  \end{eqnarray*}
%  is an equivalence of categories.  A pseudo-inverse sends $X\in \F 
% ^{I}$ to the  functor $h_{X}$ with the following description on 
% objects:
% Given an object over $K$ in $\E|I$, i.e.~a span
%   \begin{equation}\label{spann}\begin{diagram}[w=4ex,tight]
%     I & \lTo^p  & M & \rTo^q & K ,
%   \end{diagram}\end{equation}
%   then $h_{X}(p,q):= q\lowershriek(p\upperstar(X)) \in \F ^{K}$.
% \end{blanko}
\begin{proof}
We already noted that $\E | I$ is the family fibration of $\E /I 
\to \E$, which implies, cf.~(Johnstone 2002) Proposition 1.4.16 (ii), 
that it 
  is the $\E$-indexed-sum completion of $\E/I$; more precisely, if 
$\F$ has $\E$-indexed sums, then precomposition with the obvious 
``diagonal'' functor $\eta :\E /I \to \E | I$ provides an equivalence 
$$\BiFib_\E(\E|I, \F) \isopil \Fib_\E(\E/I, \F).$$ 
 We also have the well known fibered Yoneda Lemma (see for example (Borceux
 1994) Proposition 8.2.7), which says that evaluation at 
$id_{I}$ provides an equivalence 
  $$
  \Fib_\E(\E/I, \F) \isopil \F^I ,
  $$
  valid for any fibered category $\F$.   The composite of the two 
displayed functors is evaluation at $\delta$; therefore, one gets an 
explicit pseudo-inverse for the composite by composing
the pseudo-inverses for 
the two functors exhibited, and they are described in the two quoted texts (using a 
co-cleavage and a cleavage, respectively). Inspecting these 
descriptions gives in particular our (partial) description of the 
pseudo-inverse.
\end{proof}

\begin{blanko}{Corollary to the Main Lemma.}\label{cortomain}
 If $L:\E|I \to \E|J$ is a bifibered functor, then it is isomorphic, as a fibered functor, to
 $q\lowershriek \circ p\upperstar $ for the span $I \stackrel p \leftarrow M
 \stackrel q \to J$ obtained as $L(\delta)$.
\end{blanko}
\begin{proof}
  Take $\F=\E|J$ in the Main Lemma, and let $\langle p,q\rangle := L(\delta)$.
  Then we also have $q\lowershriek p\upperstar(\delta) = \langle p,q\rangle$,
  and $q\lowershriek p\upperstar$ is
  bifibered by Proposition~\ref{fibleftpreslowershriek}, since it is a fibered left adjoint.
  But the Main Lemma shows that bifibered functors $\E|I \to \E|J$
  are determined up to isomorphism by their value on $\delta$;  
  hence $L \cong q\lowershriek \circ p\upperstar $.
% 
%   Take $\F=\E|J$ in the Main Lemma, and let $\langle p,q\rangle := L(\delta)$.
%  To show that $L \cong q\lowershriek \circ p\upperstar $, the lemma implies that it
%  is enough to check at the object $\delta\in\E|I$, for which it is obvious
%  from the definition  of $p$ and $q$.
\end{proof}

\begin{cor}
 If $L:\E|I \to \E|J$ is a bifibered functor that preserves  terminal objects, then it is isomorphic, as a fibered functor,
 to $p\upperstar$, for some $p:J\to I$. \qed
\end{cor}
\begin{proof}
  The previous corollary shows that $L\cong q\lowershriek \circ 
  p\upperstar $.  Since upperstars preserve terminal objects, we 
  conclude that $q\lowershriek$ preserves the terminal object,
  but this is possible only if $q$ is invertible.  Hence
  $L \cong p\upperstar $.
\end{proof}

We record the following special case of Corollary~\ref{cortomain}:
\begin{thm}
If $L:\E|I \to \E|J$ is a fibered left adjoint, it is
obtained from a span $I \stackrel p \leftarrow M
\stackrel q \to J$, as $L\cong q\lowershriek \circ p\upperstar $.
And if $R: \E |J \to \E |I$ is a fibered right adjoint, it is 
obtained from such a span as $R\cong p\lowerstar \circ q\upperstar$.
\qed
\end{thm}

%%%%%%%%%%%%%%%%%%%%%%%%%%%%%%%%%%%%%%%%%%%%%%%%%%
\section{Local fibered right adjoints}
%%%%%%%%%%%%%%%%%%%%%%%%%%%%%%%%%%%%%%%%%%%%%%%%%%

\begin{blanko}{Local right adjoints.}
If $\F$ is a category with a terminal object $1_{\F}$, and $R:\F \to \G$ is
a functor, there is a well known canonical factorization of $R$
\begin{equation}\label{plain-Rbar}\begin{diagram}
\F & \rTo ^{\overline R} & \G /R(1_{\F}) & \rTo^{\dom} & \G
\end{diagram}\end{equation}
where $\overline R$ takes $A \in \F$ to the value of $R$ on $A\to
1_{\F}$. One says that $R$ is a {\em local right adjoint} if
$\overline{R}$ is a right adjoint. (This implies that all the evident functors
$\F /X \to \G /R(X)$ are also right adjoints, but we shall not use
this.)
\end{blanko}

\begin{blanko}{Local fibered right adjoints.}
Let $\F$ and $\G$ be fibered over $\E$, assume that $\F$ has a terminal object
$1_\F$, and let $R:\F\to\G$ be a fibered functor.  Then the factorization
\eqref{plain-Rbar} is actually a factorization of fibered functors, where
$\G/R(1_\F)$ is fibered over $\E$ via the fibration $\dom$
%%?
 and the structural functor $\G \to \E$.  If in this
situation $\overline R$ is a fibered right adjoint, we say that $R$ is a {\em
local fibered right adjoint}. (In particular, $R$ is then a (plain) local
right adjoint.)

  In the case where $\G$ is of the form $\E|J$,
  \begin{diagram}
  \F & \rTo ^{\overline R} & (\E|J) /R(1_{\F}) & \rTo^{\dom} & 
  \E|J,
  \end{diagram}
  we shall see that the middle category $(\E|J)/R(1_{\F})$ is itself 
  a fibered slice.  For this we need a little preparation:
% 
% In the case where $\G$ is of the form $\E|J$, we shall see that
% the middle object $(\E|J)/R(1_{\F})$ is itself a fibered 
% slice.  For this we need a little preparation:
\end{blanko}

\begin{blanko}{Plain slices of fibered slices.}
  In addition to the structural functor $\E|J\to\E$
  (which given a span $J \leftarrow M \to K$
  returns $K$), we have also the ``apex'' functor 
%%? 
%% (also a fibration) Do we use this; and why should it be true?
  \begin{eqnarray*}
    d: \E|J & \longrightarrow & \E  \\
    {}[J\!\leftarrow\! M \!\to\! K] & \longmapsto & M .
  \end{eqnarray*}
  For a fixed span $(J\stackrel t \leftarrow M\to K) = Q\in \E|J$,
  there is induced a forgetful fibered functor
  \begin{eqnarray*}
    (\E|J)/Q & \longrightarrow & \E|d(Q)
  \end{eqnarray*}
  which sends an object
   \begin{diagram}[w=4ex,h=3.2ex,tight, nohug]
  & & Y & \rTo & X \\
 J &\ldTo(2,1) & \dTo && \dTo \\
 & \luTo(2,1)_t & M & \rTo & K 
 \end{diagram}
 to the span 
 $$
 M \lTo Y \rTo X .
 $$
 This functor is the first leg of a factorization of the domain functor $\dom$:
 \begin{diagram}[w=4ex,h=4.5ex,tight]
   (\E|J) / Q  && \rTo^{\dom}    && \E|J    \\
 &\rdTo    &      & \ruTo_{t\lowershriek}  &  \\
 &    & \E|d(Q)    & &
 \end{diagram}
\end{blanko}

\begin{lemma}
With notation as above, when $Q$ belongs to the $1$-fiber for the structural
fibration (i.e.~is of the form $J \leftarrow M \to 1$), the forgetful functor
\begin{eqnarray*}
(\E|J)/Q & \longrightarrow & \E|d(Q)
\end{eqnarray*}
is an equivalence of fibered categories over $\E$.
\qed
\end{lemma}

We now return to the factorization \eqref{plain-Rbar} as fibered functors, and take $\F=\E|I$
and $\G=\E|J$.  Note that $\E|I$ has a terminal object $1_I$, namely
the span $I\stackrel{\id_I}\leftarrow I \to 1$. It belongs to 
the $1$-fiber, and hence so does $Q:=R(1_{I})$.  Combining the above discussion with the Lemma, we arrive at this:
\begin{cor}\label{Rt}
  Any local fibered right adjoint $R: \E|I \to \E|J$ has a canonical
factorization by fibered functors
\begin{equation}\begin{diagram} \E|I &\rTo ^{\overline{R}}& \E|B &\rTo
^{t\lowershriek}&\E|J
\end{diagram},\label{factor}\end{equation}
with $\overline R$ a fibered right adjoint.  Here
$(J \stackrel t \leftarrow B \to 1):= R^{1}(I \stackrel{\id_I}\leftarrow I \to 1)$.
\qed
\end{cor}

We can now prove the Theorem announced in the title.
% Let $\E$ be a
% locally cartesian closed category with a terminal object.
  \begin{thm}\label{mainthm}
Let $\E$ be a
locally cartesian closed category with a terminal object.
  If $R : \E|I \to \E|J$ is a local fibered right
  adjoint, then it is a polynomial functor.
  \end{thm}

  \begin{proof}
  By the description in \ref{polyparagraph} we need to construct the polynomial
% \begin{thm}\label{mainthm}
%   If $R : \E|I \to \E|J$ is a local fibered right
%    adjoint, then it is a polynomial functor (in the sense of \ref{polyparagraph}).
% \end{thm}
% 
% \begin{proof}
% We need to construct the polynomial 
  \begin{equation}\label{bridgeinproof}
\begin{diagram}[w=4ex,tight]
I& \lTo ^{s}&E& \rTo ^{p}&B&\rTo ^{t} &J 
\end{diagram}\end{equation}
representing $R$, i.e.~such that 
$R \cong t\lowershriek\circ p\lowerstar\circ s\upperstar$.
  By Corollary~\ref{Rt} we have $R = t\lowershriek \circ \overline{R}$, where
$\overline{R}: \E|I \to \E|B$ has a fibered left adjoint $L$.
Explicitly, $B$ and $t$ are determined by
$$
(J\stackrel t\leftarrow B\to 1)\ := \ R^{1}(I\stackrel{\id_I}\leftarrow I\to 1).
$$
By the Main Lemma (actually its Corollary~\ref{cortomain})
we can write $L \cong s\lowershriek \circ p\upperstar$, hence
$\overline{R}\cong p\lowerstar \circ
s\upperstar$.
% The proof of the Main Lemma gives the maps $s$ and $p$ explicitly:
The maps $s$ and $p$ are given explicitly by
$$
(I \stackrel s \leftarrow E \stackrel p \to B) \ := \ L^B( B 
\stackrel{\id_B}\leftarrow B \stackrel{\id_B}\to B) .
$$
Altogether we have $R \cong t\lowershriek\circ p\lowerstar\circ s\upperstar$
as claimed.
\end{proof}

\begin{blanko}{Remarks.}
  The converse of the theorem is also true: (fibered) polynomial functors are
  always local fibered right adjoints.  Indeed, if $P:\E|I \to \E|J$ is
  the fibered functor $P=t\lowershriek\circ p\lowerstar\circ s\upperstar $, then we have
  $\overline P = p\lowerstar \circ s\upperstar $ with fibered left adjoint
  $s\lowershriek \circ p\upperstar$.
  
  It should also be noted that the diagram \eqref{bridgeinproof} ``representing'' a
  local fibered right adjoint is essentially unique.  This follows from
  Theorem~2.17 in (Gambino-Kock 2009).  That theorem establishes a
  biequivalence between a bicategory whose $1$-cells are ``polynomials''
  \eqref{bridge} and a $2$-category whose $1$-cells are ``polynomial
  functors'' (i.e.~functors that are isomorphic to one given by a polynomial,
  and hence have a canonical extension to a fibered functor).
  Theorem~\ref{mainthm} gives an intrinsic characterization of this essential
  image.
\end{blanko}

\begin{blanko}{Example: $\E=\Set$.}
If $\E$ is the category of sets, any
category $\F$ can canonically be seen as the $1$-fiber of a category fibered
over $\E$, namely the category whose $I$-fiber is the category of
$I$-indexed families of objects in $\F$; and any functor $\F \to \G$
extends canonically to a fibered functor. One often expresses this
by saying ``any category $\F$ `is' fibered over Sets, and any
functor `is' fibered''. So rather than considering local fibered
right adjoints, one can consider just local right adjoints $\E /I \to
\E /J$ and prove that they are given by a polynomial (\ref{bridge});
see e.g.~(Gambino-Kock 2009).
\end{blanko}

If $\E$ is a more
general topos, functors $\E /I\to \E /J$ need not ``be'' fibered,
not even for $I=J=1$,
i.e.\ they may not be the $1$-fibers of a fibered functor $\E|I
\to \E|J$, as the following examples show.

\begin{blanko}{Example (Weber).}
  Let $\E$ be the category of directed graphs, i.e.~the presheaf
  category of $(0 \rightrightarrows 1)$, and let $T:\E\to\E$ be the
  free-category monad.  Weber (2007) observes that $T$ is a local right adjoint
  but not a polynomial functor.  The argument (given in detail in Example~2.5 of
  op.~cit.) amounts to showing that the left adjoint to $\overline T$
  does not preserve monos.  In contrast, for a polynomial functor
  $P=t\lowershriek p\lowerstar s\upperstar $, the left adjoint to $\overline P =
  p\lowerstar s\upperstar $ is $s\lowershriek p\upperstar $, which does preserve
  monos (as does every polynomial functor).  It follows a posteriori from our
  Main Theorem that $T$ cannot be the $1$-fiber of a local fibered right adjoint
  $\E|1 \to \E|1$.
\end{blanko}

\begin{blanko}{Example.}
  Let $\E$ be the category of $G$-sets, where $G$ is
a non-trivial group. The group homomorphism $p:G\to 1$ induces 
functors \mbox{$p\lowershriek\dashv p\upperstar \dashv p\lowerstar $,}
where $p\upperstar$ applied to a 
set $S$ makes it into a $G$-set, with trivial $G$-action. (The functors 
$p\lowershriek$ and $p\lowerstar$ may be seen as left- and right- Kan extensions, 
respectively.) The endofunctor $R: \E \to \E$ given by $p\upperstar\circ p\lowerstar$ 
converts a $G$-set $X$ to the subset consisting of the fixpoints for 
the action (and equipped with trivial action). Now $R$ has a left 
adjoint, namely $p\upperstar\circ p\lowershriek$. 
So $R$ is a right adjoint $\E \to \E$, (or $\E /1 \to \E /1$), but it is not 
the $1$-fiber of a
fibered right adjoint $\E|1 \to \E|1$. For, this would
imply, by general theory (cf.~ \ref{fib->st} below), that $R$ could be equipped with a
tensorial strength (in the sense of (Kock 1972))
$I \times R (X) \to R(I\times X)$, natural in $I$ and $X$ (objects
of $\E$). In particular, take $X=1$, so we have $I=I\times R(1) \to
R(I)$; but if $I$ is a non-empty $G$-set without stationary points,
this is impossible, since then $R(I)$ is empty.
\end{blanko}

\begin{blanko}{Example.} (Cf.~(Yetter 1987), 
(Kock-Reyes 1999) \S 4.)
  Let $D$ be an atom in a lccc $\E$, meaning that the endofunctor 
  $X\mapsto X^{D}$ has a right adjoint.  Clearly $(-)^D$
   extends to a fibered functor $\E|1 \to \E |1$ since it is 
polynomial. By (Yetter 1987), this 
extension  has a right adjoint fiber by fiber, but there cannot be a 
fibered right adjoint, unless $D=1$. 
    Indeed,  the 
  strength of $(-)^D$ is given by the obvious $I\times X^D \to (I\times X )^D$, and
  if the functor $(-)^D$ is a fibered left adjoint, this
   implies that the strength is an isomorphism 
  (cf.~\ref{fib->st}). In particular, instantiating 
  at $X=1$, one gets the ``diagonal'':
  $I \to I^D$ (exponential transpose of the projection map 
  $I\times D \to I$) which is then an isomorphism, natural 
  in $I$.  By an easy application of the strong Yoneda Lemma
  (cf.~e.g.~(Gambino-Kock 2009) Lemma 2.6), we get that $D =1$.
\end{blanko}

% \begin{blanko}{Example.} (Cf.~(Kock-Reyes 1999) \S 4.)
% Let $D$ be an atom in a lccc, meaning that
% the endofunctor $X\mapsto X^{D}$ has a right adjoint $R$. This $R$
% cannot be part of a fibered adjunction unless $D=1$ (although $R$ can
% actually be extended to a fibered functor, cf.~loc.~cit., and $L=
% (-)^{D}$ clearly
% extends to a fibered functor since it is polynomial).
% Its strength is the obvious $I\times X^D \to (I\times X )^D$; in
% particular, instantiating at $X=1$, one gets the ``diagonal'' : $I \to
% I^D$ (exponential transpose of the projection map $I\times D \to I$).
% However, the adjointness $(-)^D \dashv R$ cannot be fibered unless $D=1$; for
% (cf.~\ref{fib->st}),
% this would imply that the strength is an isomorphism,
% in particular, the diagonal $I\to I^D$ is an isomorphism, natural in
% $I$.
% By an easy application
% of the strong Yoneda Lemma (cf.~e.g.~(Gambino-Kock 2009) Lemma 2.6), we get that $D =1$.
% \end{blanko}

\begin{blanko}{From fibered functors to strength.}\label{fib->st}
  \newcommand{\proj}{\operatorname{pr}}
If $\F$ is fibered over $\E$ and has $\E$-indexed 
%%?
sums, then 
$\F^{1}$ is tensored over $\E$:  if $S$ is in $\E$ 
and $X\in \F^{1}$, we have 
$S \otimes X = \proj\lowershriek(\proj\upperstar(X ))$, where $\proj$ denotes the 
unique map $S\to 1$. If $L: \F \to \G$ is a fibered functor
(where $\G$ likewise has $\E$-indexed 
%%?
sums), we may 
rewrite $S\otimes L^{1}(X )$ as 
$$
S\otimes L^{1}(X ) = \proj\lowershriek(\proj\upperstar(L^{1}(X)))
\cong\proj\lowershriek(L^{S} (\proj\upperstar(X))).
$$
On the other hand, the $2$-cell exhibited in (\ref{lowershriek2cell}) in
particular provides a map 
$$\proj\lowershriek (L^{S}(\proj\upperstar  X ) 
\to L^{1} (\proj\lowershriek\proj\upperstar X )=L^{1}(S\otimes X ).$$
 So by composition, 
we get a map
$$
S\otimes L^{1}(X ) \to L^{1}(S\otimes X),
$$
which is a tensorial strength for the functor $L^{1}$.   So briefly, 
``a fibering
implies a (tensorial) strength'', cf.~(Johnstone 1997)~\S3 for the case
$\F = \G = \E$.  Note that if $L$ is a fibered left adjoint, it commutes with
lowershriek (cobase-change) functors, and therefore the strength is an isomorphism.
\end{blanko}

\begin{blanko}{Is strength sufficient?}
There is a partial converse to ``fibering implies strength'', due to Par\'{e},
cf.~(Johnstone 1997) Proposition 3.3: {\em a pullback-preserving functor
$\E\to\E$ with a strength extends uniquely to a fibered functor $\E|1 \to \E|1$}.
It is plausible that a similar result could hold also for
functors of the form $\E/I \to \E/J$, and that in this situation strong
natural transformation would correspond to fibered natural
transformations, so that the notions of strong adjunction and
fibered adjunction would match up.

In (Gambino-Kock 2009), polynomial functors are studied in the
setting of functors equipped with a tensorial strength (``strong
functors'') and strong natural transformations; in a sense this
is more economical, compared to the fibered setting, since only
plain slices are needed.  Should the above speculation turn true,
it would seem
to imply a strength version of our main theorem, namely
that every local strong right adjoint $\E/I \to \E/J$ is
polynomial.
\end{blanko}

\noindent
{\bf Acknowledgments.}
  The authors thank Martin Hyland for insightful and pertinent comments.
% %%?
%   The second-named author was partially supported by research grants 
%   MTM2009-10359 % Nart
%        and
% MTM2010-20692 % Castellana
% of the Spanish Ministry of Science and Innovation.

  \makeatletter
  
%   \newdimen\bibindent
% \setlength\bibindent{1.5em}
% \renewenvironment{thebibliography}[1]
%      {\section*{\refname}%
% %       \@mkboth{\MakeUppercase\refname}{\MakeUppercase\refname}%
%       \list{}%\@biblabel{\@arabic\c@enumiv}}%
%            {\settowidth\labelwidth{0em}%
% %             \leftmargin\labelwidth
% %             \advance\leftmargin\labelsep
% %             \@openbib@code
%             \usecounter{enumiv}%
% %             \let\p@enumiv\@empty
% %             \renewcommand\theenumiv{\@arabic\c@enumiv}}%
%     }  \sloppy
%       \clubpenalty4000
%       \@clubpenalty \clubpenalty
%       \widowpenalty4000%
%       \sfcode`\.\@m}
%      {\def\@noitemerr
%        {\@latex@warning{Empty `thebibliography' environment}}%
%       \endlist}

\renewenvironment{thebibliography}[1]
 {\section*{\refname}%
   \addcontentsline{toc}{section}{\refname}%
   \list{}{\labelwidth\z@ \leftmargin 1em \itemindent-\leftmargin}%
   \small\normalfont \parindent\z@
   \parskip\z@ \@plus .1\p@\relax
   \sloppy\clubpenalty\z@ \widowpenalty\@M
   \sfcode`\.\@m\relax}
  {\def\@noitemerr
       {\@latex@warning{Empty `thebibliography' environment}}%
   \endlist}

      \makeatother

\medskip

\noindent
{\sc Anders Kock} \url{<kock@imf.au.dk>}
\\
Matematisk Institut\\
Aarhus Universitet\\
Denmark

\bigskip

\noindent
{\sc Joachim Kock} \url{<kock@mat.uab.cat>}
\\
Departament de Matem\`atiques\\
Universitat Aut\`onoma de Barcelona\\
Spain

\end{document}